\newcommand{\me}{\mathrm{e}}
\newcommand{\mi}{\mathrm{i}}
\newcommand{\dif}{\mathrm{d}}
\newcommand{\Dif}{\mathrm{D}}

\newcommand{\cir}{_{\mathrm{c}}}
\newcommand{\ave}{_{\mathrm{ave}}}
\newcommand{\psiA}{\psi_{\mathrm{A}}}
\newcommand{\psiL}{\psi_{\mathrm{L}}}
\newcommand{\psiK}{\psi_{\mathrm{K}}}
\newcommand{\psiAO}{\psi_{\mathrm{A0}}}
\newcommand{\tr}{^\mathrm{T}}
\newcommand{\R}{\mathrm{Re}}
\newcommand{\I}{\mathrm{Im}}
\newcommand\p{\ensuremath{\partial}}

\documentclass[preprint,12pt]{elsarticle}
\let\today\relax
\makeatletter
\def\ps@pprintTitle{%
    \let\@oddhead\@empty
    \let\@evenhead\@empty
    \def\@oddfoot{\footnotesize\itshape
         {Accepted by SCIENTIA SINICA Mathematica (in Chinese) } \hfill\today}%
    \let\@evenfoot\@oddfoot
    }
\makeatother
\usepackage{amsthm}
\usepackage{amssymb}
\usepackage[tbtags,fleqn]{amsmath}
\usepackage{bm}
\usepackage[hypertex]{hyperref}
\usepackage{floatrow}
\usepackage{graphicx}
\usepackage[label font=bf,labelformat=simple]{subfig}
\usepackage{caption}

\newtheorem{thm}{Theorem}[section]

\newtheorem{lem}[thm]{Lemma}
\newdefinition{rmk}[thm]{Remark}

\begin{document}
\date{Febrauray 04,2021}
\begin{frontmatter}



\title{Energy stability of the Charney--DeVore quasi-geostrophic equation  for atmospheric blocking}

\author[a1]{Zhi-Min Chen}
\ead{zmchen@szu.edu.cn}
\address[a1]{School  of Mathematics and Statistics, Shenzhen University, Shenzhen 518060,  China}%

\author[a2]{Xiangming Xiong\corref{cor}}
\ead{xxiong@princeton.edu}
\address[a2]{The Program in Applied and Computational Mathematics, Princeton University,  Princeton, NJ 08544, USA}

\cortext[cor]{Corresponding author.}

\begin{abstract}
Charney and DeVore [J. Atmos. Sci. 36 (1979), 1205-1216] found multiple equilibrium states as a consequence of bottom topography in their pioneering work on the quasi-geostrophic barotropic flow over topography in a $\beta$-plane channel. In the present paper, we prove that the basic flow is asymptotically stable in a parameter region, including the flat topography situation, which excludes the existence of multiple equilibrium states therein. Moreover, we show that an additional condition on the average zonal force or the average zonal velocity is indispensable to the well-posedness of the Charney--DeVore quasi-geostrophic equation. Coexistence of at least three equilibrium states is confirmed by a pseudo-arclength continuation method for different topographic amplitudes. The stabilities of the equilibrium states are examined by high-resolution direct numerical simulations.
\end{abstract}


\end{frontmatter}


\section{Introduction}
\label{sec:introduction}

A surface quasi-geostrophic flow is the first order approximation of a geostrophic flow with
respect to the quasi-geostrophic approximation under a small Rossby number \cite{Pedlosky1987}. From the mathematical formulation of \cite{ConstantinMajdaTabak1994}, the energy dissipation
controlled by the fractional Laplacian $\kappa (-\Delta)^\alpha$ is applied to a quasi-geostrophic
flow \cite{Wu1997,ConstantinWu1999}. Mathematical theory of  the dissipative quasi-geostrophic
equation has been extensively studied (see, for example, \cite{ChenMiaoZhang2007,Kiselev2007,Wu1997,ConstantinIyerWu2008,Wu2005,ChenPrice2008,Chen2016}). In the present study, we are interested in a quasi-geostrophic flow, which is from the understanding of blocks in atmosphere and involves  the weak energy dissipation $\kappa (-\Delta)^\alpha$ with $\alpha=0$.

Large-scale flows in the atmosphere are strongly influenced by the bottom topography. When the topography is not flat, the small-scale and large-scale components of atmospheric flow interact via not only advection but also topographic stress. Charney and DeVore \cite{Charney1979} proposed a barotropic channel model to study the influence of topography on the atmospheric flow. They found multiple equilibrium states in the quasi-geostrophic $\beta$-plane flow over topography with highly truncated spectral expansions. In some region of parameter space, there were one unstable equilibrium state and two stable equilibrium states. The two stable equilibrium states comprised a ``high-index'' flow with a relatively larger zonal flux and a ``low-index'' flow showing atmospheric blocking. Atmospheric blocking refers to the unusual persistence of a large-scale weather pattern in an area, and usually occurs when the normal midlatitude eastward flow is interrupted by topography or strong meridional flow \cite{Ghil1987, Cushman2011}. Atmospheric blocking may bring about extreme weather events such as floods, droughts, and persistent abnormal temperatures \cite{Dolzhansky2013}.

In the work of Charney and DeVore \cite{Charney1979}, the governing equation of the vertical vorticity $\zeta$ was
\begin{equation}
\frac{\p \zeta}{\p t} + \frac{\p \psi}{\p x} \frac{\p (\zeta + h + \beta y)}{\p y} - \frac{\p \psi}{\p y} \frac{\p (\zeta + h + \beta y)}{\p x} = -k \zeta - \frac{\p F}{\p y},\label{eq:psi_CD}
\end{equation}
where $\psi$ was the stream function satisfying $\zeta=\p^2 \psi/\p x^2 + \p^2 \psi/\p y^2$, and $h$ was the bottom topography. $\beta=2\Omega(\cos\phi_0)/a$, where $\Omega$ was the angular speed of the earth's rotation, $a$ was the radius of the earth, and $\phi_0$ was the average latitude of the zonal channel. $k$ was a coefficient of friction exerted by the Ekman layer near the ground.  $-\p F/\p y$ was the vorticity source arisen from the zonal body force $F$.

Charney and DeVore \cite{Charney1979} assumed the flow was periodic in the zonal direction with a non-dimensional period of $2\pi$. The meridional boundaries $y=0$ and $y=\pi$ were assumed to be non-penetrative. Therefore, the boundary conditions of the stream function $\psi$ were
\begin{equation}
\psi(x+2\pi,y,t)  = \left.\left.\psi(x,y,t) ,\qquad \frac{\p \psi}{\p x}\right |_{y=0} = \frac{\p \psi}{\p x} \right|_{y=\pi} = 0.\label{eq:bc_CD}
\end{equation}

However, the governing equation (\ref{eq:psi_CD}) was not well-posed under the boundary conditions (\ref{eq:bc_CD}). Although the average zonal force did not appear in the governing equation (\ref{eq:psi_CD}), it evidently influenced the atmospheric flow. Therefore, another equation of the average zonal velocity is needed to take into account this influence. This problem was noticed by Davey \cite{Davey1980}. He proposed a quasi-linear theory for steady flow over topography in a periodic channel, which included the equation of the average zonal velocity but neglected the advection due to the mean shear and the interactions between eddies. He confirmed that there were multiple equilibrium states when the amplitude of topography was sufficiently large.

The multiple equilibrium states of the Charney--DeVore quasi-geostrophic equation were analysed and calculated by many researchers. Hart \cite{Hart1979} found stable equilibrium states when the topography was anisotropic and mainly varied in the zonal direction. In the fully-nonlinear calculation performed by Legras and Ghil \cite{Legras1985}, however, the parameter region where multiple equilibrium states existed was found to be smaller than that in a quasi-linear model. Tung and Rosenthal found that the multiple equilibrium states only existed in a parameter region that was far from realistic if fully-nonlinear effects were considered and the truncation of modes was not severe \cite{Tung1985}. However, multiple equilibrium states were still found with wind speeds only twice as large as the observed ones by introducing mid-latitudinal zonal jets in the model \cite{Zidikheri2007}. Other researchers focused on the more complicated case with the baroclinic model \cite{Linden1983} or two layer flow \cite{Vallis1985}.  Multiple equilibrium states from  a bifurcation  viewpoint of a non-parallel vortex flow were investigated in \cite{Chen2019}. One may also refer to \cite{Chen2002} for the Hopf bifurcation of a quasi-geostrophic   atmospheric  flow involving strong energy dissipation.

When the potential vorticity is proportional to the stream function, the nonlinear Jacobian term in the equation of potential vorticity will vanish. The inviscid stabilities of the equilibrium states with potential vorticity proportional to the stream function have already been analysed in a $\beta$-plane channel \cite{Fyfe1986, Fyfe1989}. In the absence of forcing and dissipation, the nonlinear stability of the equilibrium states was studied with the conservation of a combinition of the kinetic energy and the potential enstrophy in a double-periodic domain \cite{Carnevale1987, Frederiksen2005}.
Zou and Fyfe considered the energy stability of the equilibrium states of forced and dissipated barotropic flow \cite{Zou1993}. Energy stability theory of general parallel flows was recently investigated by the authors \cite{XiongChen2019}.

In \cite{Chen2016b}, we proved the global stability of the equilibrium state in a parameter region, and performed numerical simulations with 12 Fourier modes in the zonal direction and 65 collocation points in the meridional direction. In this paper, we will prove the global stability of the equilibrium state in a larger parameter region. After that, we will calculate the equilibrium states with a pseudo-arclength continuation method, and will examine their stabilities with direct numerical simulations with higher resolution. The derivation of the governing equations of the quasi-geostrophic barotropic flow is presented in Section \ref{sec:model}. The proof of the global stability is given in Section \ref{sec:theorem}. After briefly introducing the numerical methods in Section \ref{sec:method}, we present the numerical results of multiple equilibrium states in Section \ref{sec:results}. Discussions on the equilibrium states found by Charney and DeVore are made in Section \ref{sec:discussion}, and conclusions are given in Section \ref{sec:conclusion}.

\section{Mathematical model}
\label{sec:model}

The large-scale flow of atmosphere outside the Ekman boundary layer is nearly inviscid, and therefore the governing equations and boundary conditions are
\begin{align}
& \nabla^* \cdot \bm{U}^*=0,\\
& \frac{\p \bm{U}^*}{\p t^*} + (\bm{U}^* \cdot \nabla^*)\bm{U}^* + \bm{f}^* \times \bm{U}^* = -\frac{1}{\rho} \nabla^* P^* + \frac{1}{\rho} \bm{F}^*,\\
& (\bm{U}^*,P^*) (\bm{x}^*+2\pi L\bm{e}_x,t^*) = (\bm{U}^*,P^*) (\bm{x}^*,t^*),\\
&\left. V^*\right|_{y^*=0} =\left. V^*\right|_{y^*=\pi L} = 0,\\
& W^*|_{z^*=h^*} = U^*|_{z^*=h^*} \frac{\p h^*}{\p x^*} + V^*|_{z^*=h^*} \frac{\p h^*}{\p y^*},\ W^*|_{z^*=H} = 0, \label{eq:bc_z}
\end{align}
where $\bm{x}^*=(x^*,y^*,z^*)$, and $x^*$ and $y^*$ are the zonal (eastward) and meridional (northward) coordinates in the tangent plane of the earth, respectively. The coordinate in the direction perpendicular to the tangent plane (opposite to the direction of the gravitational acceleration) is $z^*$. The unit vectors in these directions are $\bm{e}_x$, $\bm{e}_y$, and $\bm{e}_z$. Assume the latitude $\phi$ to be confined to a small interval around a middle latitude $\phi_0$ in the Northern Hemisphere, and then we have $x^* \approx a\theta \cos \phi_0$ and $y^* \approx \pi L/2 + a(\phi-\phi_0)$, where $\theta$, $a$, and $\pi L$ are the longitude, the radius of the earth, and the distance between zonal walls, respectively. $\bm{U}^*=(U^*,V^*,W^*)$ is the velocity of the atmosphere. The density of the atmosphere $\rho$ is assumed to be a constant. $H$ is the height of the upper free surface, and $h^*(x^*,y^*)$ is the lower topography, which satisfies
\begin{equation*}\iint h^* \ \dif x^*\dif y^*=0.\end{equation*}

According to the $\beta$-plane approximation, $\bm{f}^* \approx 2\Omega (\sin\phi) \bm{e}_z \approx (f^*_0 + \beta^* y^*)\bm{e}_z$, where $\Omega$ is the angular speed of the earth's rotation, $f^*_0 = 2\Omega\sin\phi_0 - \Omega\pi L(\cos\phi_0)/a$, and $\beta^*=2\Omega(\cos\phi_0)/a$. Only the local normal component of the earth's rotation is preserved here \cite{Pedlosky1987}. The body force is $\bm{F}^* = F^* \bm{e}_x - \rho g \bm{e}_z$, where $g$ is the magnitude of the gravitational acceleration, and $F^*$ is the external force acting in the zonal direction. In the above definitions, asterisks denote dimensional variables. The constants $H$, $L$, $a$, $\Omega$, $\rho$, and $g$ are also dimensional quantities.

Suppose $|h^*| \ll H \ll \pi L \ll 2a$, and therefore the vertical scale of the geophysical flow is far less than its horizontal scale. Neglecting $W^*$, $\p U^*/\p z^*$, and $\p V^*/\p z^*$ in the momentum equations, we have
\begin{align}
& \frac{\p U^*}{\p x^*} + \frac{\p V^*}{\p y^*} + \frac{\p W^*}{\p z^*} =0,\label{eq:continuity}\\
& \frac{\p U^*}{\p t^*} + U^* \frac{\p U^*}{\p x^*} + V^* \frac{\p U^*}{\p y^*} - (f^*_0 + \beta^* y^*)V^* = -\frac{1}{\rho} \frac{\p P^*}{\p x^*} + \frac{1}{\rho} F^*,\label{eq:momentum_x}\\
& \frac{\p V^*}{\p t^*} + U^* \frac{\p V^*}{\p x^*} + V^* \frac{\p V^*}{\p y^*} + (f^*_0 + \beta^* y^*)U^* = -\frac{1}{\rho} \frac{\p P^*}{\p y^*},\label{eq:momentum_y}\\
& 0= -\frac{1}{\rho} \frac{\p P^*}{\p z^*} - g. \label{eq:momentum_z}
\end{align}
The vertical momentum equation (\ref{eq:momentum_z}) implies
\begin{equation*}
\frac{\p}{\p z^*}(\frac{\p P^*}{\p x^*}) = \frac{\p}{\p z^*}(\frac{\p P^*}{\p y^*}) = 0.
\end{equation*}
Assume $F^*= F^*(x^*, y^*)$, and then all terms in the horizontal momentum equations (\ref{eq:momentum_x})--(\ref{eq:momentum_y}) are independent of $z^*$.

Integrating the continuity equation (\ref{eq:continuity}) from $z^*=h^*$ to $z^*=H$, we have
\begin{equation*}
W^*|_{z^*=H} - W^*|_{z^*=h^*} = -\int^H_{h^*} (\frac{\p U^*}{\p x^*} + \frac{\p V^*}{\p y^*})\ \dif z^* = -(H-h^*)(\frac{\p U^*}{\p x^*} + \frac{\p V^*}{\p y^*}).
\end{equation*}
By substituting the boundary condition (\ref{eq:bc_z}) into the above equation, we obtain
\begin{equation*}
\frac{\p}{\p x^*} ((H-h^*)U^*) + \frac{\p}{\p y^*} ((H-h^*)V^*) = 0.
\end{equation*}

Let $\tilde{U}^* = (1-h^*/H) U^*$ and $\tilde{V}^* = (1-h^*/H) V^*$, and then the equations (\ref{eq:continuity})--(\ref{eq:momentum_y}) are approximated by
\begin{align}
& \frac{\p \tilde{U}^*}{\p x^*} + \frac{\p \tilde{V}^*}{\p y^*} = 0,\label{eq:continuity_tilde}\\
& (1+\frac{h^*}{H})\Big( \frac{\p \tilde{U}^*}{\p t^*} + (\tilde{U}^* \frac{\p}{\p x^*} + V^* \frac{\p}{\p y^*})(\tilde{U}^* +\frac{h^*}{H} \tilde{U}^*) - (f^*_0 + \beta^* y^*)\tilde{V}^* \Big) = -\frac{1}{\rho} \frac{\p P^*}{\p x^*} + \frac{1}{\rho} F^*,\label{eq:momentum_x_tilde}\\
& (1+\frac{h^*}{H})\Big( \frac{\p \tilde{V}^*}{\p t^*} + (\tilde{U}^* \frac{\p}{\p x^*} + V^* \frac{\p}{\p y^*})(\tilde{V}^* +\frac{h^*}{H} \tilde{V}^*) + (f^*_0 + \beta^* y^*)\tilde{U}^* \Big) = -\frac{1}{\rho} \frac{\p P^*}{\p y^*},\label{eq:momentum_y_tilde}
\end{align}
where $(1-h^*/H)^{-1} \approx 1+h^*/H$ is used.

The characteristic length, velocity, and time of large-scale atmospheric motions are $L \sim 10^3\, \mathrm{km}$, $U^*_0 \sim 10\, \mathrm{m/s}$, and $L/U^*_0 \sim 10^5\, \mathrm{s} \approx 28\, \mathrm{hours}$, respectively \cite{Vallis2006, Holton2013}. The average depth of the atmosphere and the characteristic elevation of topography are supposed to be $H \sim 10\, \mathrm{km}$ \cite{Holton2013} and $h^*_0 \sim 1\, \mathrm{km}$, respectively. Suppose $\phi_0 \approx 45^\circ \, \mathrm{N}$. Using $a \approx 6400\, \mathrm{km}$, $g \approx 10\, \mathrm{m/s^2}$, and $\Omega \approx 2\pi/86400\, \mathrm{s^{-1}}$, we have $f^*_0 \sim 10^{-4}\, \mathrm{s^{-1}}$ and $\beta^* \sim 10^{-11}\, \mathrm{m^{-1}s^{-1}}$. The Rossby number is $Ro = U^*_0/(f^*_0 L) \sim 0.1$. Then we have the following estimates:
\begin{align*}
& f^*_0 \tilde{U}^*, f^*_0 \tilde{V}^* \sim f^*_0 U^*_0 \sim 10^{-3}\, \mathrm{m/s^2},\\
& \frac{h^*}{H} f^*_0 \tilde{U}^*, \frac{h^*}{H} f^*_0 \tilde{V}^* \sim \frac{h^*}{H} f^*_0 U^*_0 \sim 10^{-4}\, \mathrm{m/s^2},\\
& \beta^* y^* \tilde{U}^*, \beta^* y^* \tilde{V}^* \sim \beta^* L U^*_0 \sim 10^{-4}\, \mathrm{m/s^2},\\
& \frac{\p\tilde{U}^*}{\p t^*}, \frac{\p\tilde{V}^*}{\p t^*} \sim \frac{U^{*2}_0}{L} \sim 10^{-4}\, \mathrm{m/s^2},\\
& \tilde{U}^* \frac{\p\tilde{U}^*}{\p x^*}, \tilde{V}^* \frac{\p\tilde{U}^*}{\p y^*}, \tilde{U}^* \frac{\p\tilde{V}^*}{\p x^*},  \tilde{V}^* \frac{\p\tilde{V}^*}{\p y^*} \sim \frac{U^{*2}_0}{L} \sim 10^{-4}\, \mathrm{m/s^2}.
\end{align*}
It follows that the terms on the right hand side of the equations (\ref{eq:momentum_x_tilde})--(\ref{eq:momentum_y_tilde}) are on the order of $10^{-3}\, \mathrm{m/s^2}$.

Neglecting terms that have magnitudes less than $10^{-4}\, \mathrm{m/s^2}$ in the equations (\ref{eq:continuity_tilde})--(\ref{eq:momentum_y_tilde}), we have
\begin{align}
& \frac{\p \tilde{U}^*}{\p x^*} + \frac{\p \tilde{V}^*}{\p y^*} = 0,\label{eq:continuity_tilde_without_friction}\\
& \frac{\p \tilde{U}^*}{\p t^*} + \tilde{U}^* \frac{\p \tilde{U}^*}{\p x^*} + V^* \frac{\p \tilde{U}^*}{\p y^*} - (f^*_0 + \frac{f^*_0 h^*}{H} + \beta^* y^*)\tilde{V}^* = -\frac{1}{\rho} \frac{\p P^*}{\p x^*} + \frac{1}{\rho} F^*,\\
& \frac{\p \tilde{V}^*}{\p t^*} + \tilde{U}^* \frac{\p \tilde{V}^*}{\p x^*} + V^* \frac{\p \tilde{V}^*}{\p y^*} + (f^*_0 + \frac{f^*_0 h^*}{H} + \beta^* y^*)\tilde{U}^* = -\frac{1}{\rho} \frac{\p P^*}{\p y^*}.\label{eq:momentum_y_tilde_without_friction}
\end{align}

In the above derivation, we focus on the large-scale flow of atmosphere outside the Ekman boundary layer. However, the ground exerts friction to the flow of atmosphere through the Ekman boundary layer. Denote the zonal velocity, the meridional velocity, and the pressure in the Ekman boundary layer as $\tilde{U}^*_E(z^*)$, $\tilde{V}^*_E(z^*)$, and $P^*_E(x^*, y^*, z^*)$, respectively. The governing equations and boundary conditions are
\begin{align}
& - f^*_0 \tilde{V}^*_E = -\frac{1}{\rho} \frac{\p P^*_E}{\p x^*} + \frac{1}{\rho} F^* + \nu_E \frac{\p^2 \tilde{U}^*_E}{\p z^{*2}},\label{eq:momentum_x_Ekman}\\
& f^*_0 \tilde{U}^*_E = -\frac{1}{\rho} \frac{\p P^*_E}{\p y^*} + \nu_E \frac{\p^2 \tilde{V}^*_E}{\p z^{*2}},\label{eq:momentum_y_Ekman}\\
& 0 = -\frac{1}{\rho} \frac{\p P^*_E}{\p z^*} - g,\label{eq:momentum_z_Ekman}\\
& \tilde{U}^*_E = \tilde{V}^*_E = 0, \quad z^*=0,\\
& \tilde{U}^*_E \rightarrow \tilde{U}^*, \tilde{V}^*_E \rightarrow \tilde{V}^*, P^*_E \rightarrow P^*, \quad \delta^*_E \ll z^* \ll H,
\end{align}
where topography is neglected for simplicity. $\delta^*_E$ is the characteristic depth of the boundary layer, which is to be determined. $\nu_E$ is the bulk eddy viscosity in the boundary layer. From the equation (\ref{eq:momentum_z_Ekman}), we have
\begin{equation*}
\frac{\p}{\p z^*}(\frac{\p P^*_E}{\p x^*}) = \frac{\p}{\p z^*}(\frac{\p P^*_E}{\p y^*}) = 0.
\end{equation*}
As a result,
\begin{align}
& -\frac{1}{\rho}\frac{\p P^*_E}{\p x^*} + \frac{1}{\rho} F^* = -\frac{1}{\rho}\frac{\p P^*_E}{\p x^*} |_{z^* \gg \delta^*_E} + \frac{1}{\rho} F^* = -\frac{1}{\rho}\frac{\p P^*}{\p x^*} + \frac{1}{\rho} F^* \approx -f^*_0 \tilde{V}^*,\label{eq:dPdx}\\
& -\frac{1}{\rho}\frac{\p P^*_E}{\p y^*} = -\frac{1}{\rho}\frac{\p P^*_E}{\p y^*} |_{z^* \gg \delta^*_E} = -\frac{1}{\rho}\frac{\p P^*}{\p y^*} \approx f^*_0 \tilde{U}^*,\label{eq:dPdy}
\end{align}
where the geostrophic balance approximation is used. The geostrophic balance approximation results from the balance of terms that have magnitude of $10^{-3}\, \mathrm{m/s^2}$ in the equations (\ref{eq:continuity_tilde})--(\ref{eq:momentum_y_tilde}).

Substituting (\ref{eq:dPdx})--(\ref{eq:dPdy}) to the equations (\ref{eq:momentum_x_Ekman})--(\ref{eq:momentum_y_Ekman}), we have
\begin{align}
& -f^*_0 (\tilde{V}^*_E - \tilde{V}^*) \approx \nu_E \frac{\p^2 \tilde{U}^*_E}{\p z^{*2}},\\
& f^*_0 (\tilde{U}^*_E - \tilde{U}^*) \approx \nu_E \frac{\p^2 \tilde{V}^*_E}{\p z^{*2}}.
\end{align}
The solution is
\begin{align}
& \tilde{U}^*_E \approx \tilde{U}^* - \tilde{U}^* \exp(-\frac{z^*}{\delta^*_E})\cos(\frac{z^*}{\delta^*_E}) - \tilde{V}^* \exp(-\frac{z^*}{\delta^*_E})\sin(\frac{z^*}{\delta^*_E}),\\
& \tilde{V}^*_E \approx \tilde{V}^* - \tilde{V}^* \exp(-\frac{z^*}{\delta^*_E})\cos(\frac{z^*}{\delta^*_E}) + \tilde{U}^* \exp(-\frac{z^*}{\delta^*_E})\sin(\frac{z^*}{\delta^*_E}),
\end{align}
where $\delta^*_E = (2\nu_E/ f^*_0)^{1/2}$ \cite{Vallis2006}. Therefore, the shear stresses exerted by the ground are
\begin{align}
& \tau^*_x = -\rho \nu_E \frac{\p \tilde{U}^*_E}{\p z^*}|_{z^*=0} \approx -\frac{\rho\nu_E(\tilde{U}^* - \tilde{V}^*)}{\delta^*_E},\\
& \tau^*_y = -\rho \nu_E \frac{\p \tilde{V}^*_E}{\p z^*}|_{z^*=0} \approx -\frac{\rho\nu_E(\tilde{U}^* + \tilde{V}^*)}{\delta^*_E}.
\end{align}

The original equations (\ref{eq:continuity_tilde_without_friction})--(\ref{eq:momentum_y_tilde_without_friction}) are modified as
\begin{align}
& \frac{\p \tilde{U}^*}{\p x^*} + \frac{\p \tilde{V}^*}{\p y^*} = 0,\label{eq:continuity_tilde_with_friction}\\
& \frac{\p \tilde{U}^*}{\p t^*} + \tilde{U}^* \frac{\p \tilde{U}^*}{\p x^*} + V^* \frac{\p \tilde{U}^*}{\p y^*} - (f^*_0 + \frac{f^*_0 h^*}{H} + \beta^* y^*)\tilde{V}^* = -\frac{1}{\rho} \frac{\p P^*}{\p x^*} -\frac{\nu_E(\tilde{U}^* - \tilde{V}^*)}{\delta^*_E H} + \frac{1}{\rho} F^*,\\
& \frac{\p \tilde{V}^*}{\p t^*} + \tilde{U}^* \frac{\p \tilde{V}^*}{\p x^*} + V^* \frac{\p \tilde{V}^*}{\p y^*} + (f^*_0 + \frac{f^*_0 h^*}{H} + \beta^* y^*)\tilde{U}^* = -\frac{1}{\rho} \frac{\p P^*}{\p y^*} -\frac{\nu_E(\tilde{U}^* + \tilde{V}^*)}{\delta^*_E H},\label{eq:momentum_y_tilde_with_friction}
\end{align}
so that the volume integrals in the computational domain of the body forces $\tau^*_x/H$ and $\tau^*_y/H$ are equal to the surface integrals of the shear stresses on the ground.

Introducing $\tilde{P}^* = P^* - \rho ( f^*_0 + \nu_E/(\delta^*_E H) )\tilde{\psi}^*$, where $\tilde{\psi}^*$ is the stream function satisfying $\tilde{U}^* = -\p \tilde{\psi}^*/\p y^*$ and $\tilde{V}^* = \p \tilde{\psi}^*/\p x^*$, we have
\begin{align}
& \frac{\p \tilde{U}^*}{\p x^*} + \frac{\p \tilde{V}^*}{\p y^*} = 0,\label{eq:continuity_tilde_final}\\
& \frac{\p \tilde{U}^*}{\p t^*} + \tilde{U}^* \frac{\p \tilde{U}^*}{\p x^*} + V^* \frac{\p \tilde{U}^*}{\p y^*} - (\frac{f^*_0 h^*}{H} + \beta^* y^*)\tilde{V}^* = -\frac{1}{\rho} \frac{\p \tilde{P}^*}{\p x^*} -\frac{\nu_E}{\delta^*_E H}\tilde{U}^* + \frac{1}{\rho} F^*,\\
& \frac{\p \tilde{V}^*}{\p t^*} + \tilde{U}^* \frac{\p \tilde{V}^*}{\p x^*} + V^* \frac{\p \tilde{V}^*}{\p y^*} + (\frac{f^*_0 h^*}{H} + \beta^* y^*)\tilde{U}^* = -\frac{1}{\rho} \frac{\p \tilde{P}^*}{\p y^*} -\frac{\nu_E}{\delta^*_E H}\tilde{V}^*.\label{eq:momentum_y_tilde_final}
\end{align}

Define the non-dimensional coordinates, time, velocity, stream function, pressure, topography, and zonal force as in \cite{Charney1979}:
\begin{equation*}
(x,y) = \frac{(x^*,y^*)}{L}, \qquad t = f^*_0 t^*, \qquad (U,V)=\frac{(\tilde{U}^*,\tilde{V}^*)}{f^*_0 L}, \qquad \psi = \frac{\tilde{\psi}^*}{f^*_0 L^2}, \qquad P=\frac{\tilde{P}^*}{\rho f^{*2}_0 L^2},
\end{equation*}
\begin{equation*}
h = \frac{h^*}{H}, \qquad F = \frac{F^*}{\rho f^{*2}_0 L}.
\end{equation*}
Define non-dimensional parameters $\beta$ and $k$ as
\begin{equation*}
\beta = \frac{L}{f^*_0} \beta^* = \frac{2L}{2a\tan\phi_0 - \pi L},\qquad k = \frac{\delta^*_E}{2H}.
\end{equation*}
The zonal force $F$ is further decomposed as the sum of the average zonal force $F\ave$ and the fluctuating zonal force $F'$, which are defined as
\begin{equation*}
F\ave = \frac{1}{2\pi^2} \int_0^{\pi} \int_0^{2\pi} F \ \dif x \dif y
\end{equation*}
and $F' = F - F\ave$.

Then we have the non-dimensional governing equations and boundary conditions of the primitive variables:
\begin{align}
& \frac{\p U}{\p x} + \frac{\p V}{\p y} = 0,\label{eq:continuity_final}\\
& \frac{\p U}{\p t} + U \frac{\p U}{\p x} + V \frac{\p U}{\p y} - (h + \beta y)V = -\frac{\p P}{\p x} -kU + F\ave + F',\label{eq:momentum_x_final}\\
& \frac{\p V}{\p t} + U \frac{\p V}{\p x} + V \frac{\p V}{\p y} + (h + \beta y)U = -\frac{\p P}{\p y} -kV,\label{eq:momentum_y_final}\\
& (U,V,P)(x+2\pi,y,t) = (U,V,P)(x,y,t), \qquad V|_{y=0} = V|_{y=\pi} = 0.\label{eq:bc_final}
\end{align}

Introducing the vertical component of the non-dimensional vorticity $\zeta = \nabla^2 \psi$, where $\nabla^2 = \p^2/\p x^2 + \p^2/\p y^2$, we have the governing equation of the vertical vorticity \cite{Charney1979}:
\begin{equation}
\frac{\p \zeta}{\p t} + \frac{\p \psi}{\p x} \frac{\p (\zeta + h + \beta y)}{\p y} - \frac{\p \psi}{\p y} \frac{\p (\zeta + h + \beta y)}{\p x} = -k \zeta - \frac{\p F'}{\p y},\label{eq:psi_final}
\end{equation}
and the boundary conditions of the stream function \cite{Charney1979}:
\begin{equation}
\psi(x+2\pi,y,t)  = \psi(x,y,t) ,\left.\left.\qquad \frac{\p \psi}{\p x} \right|_{y=0} = \frac{\p \psi}{\p x} \right|_{y=\pi} = 0.\label{eq:psi_bc_final}
\end{equation}

The equations of the primitive variables (\ref{eq:continuity_final})--(\ref{eq:bc_final}), however, cannot be recovered from the equations of the vertical vorticity (\ref{eq:psi_final})--(\ref{eq:psi_bc_final}) without the information of the average zonal force $F\ave$. Therefore, the equations of the vertical vorticity need an additional condition on the average zonal force to be well-posed.

This additional condition can also be given in the form  with respect to  the average  velocity
\begin{equation}
U\ave = \frac{1}{2\pi^2} \int_0^{\pi} \int_0^{2\pi} U \ \dif x \dif y = \frac{1}{\pi} \int_0^{\pi} U \ \dif y.
\end{equation}
Observing that  \begin{equation*}\int^{2\pi}_{0} Vdx =\int^{2\pi}_{0} \frac{\partial \psi}{\partial x}dx =0
\end{equation*} due to the periodic boundary condition \eqref{eq:psi_bc_final},  we notice from the equation (\ref{eq:momentum_x_final}) that the average zonal velocity and the average zonal force are coupled through
\begin{equation}
\frac{\dif U\ave}{\dif t} = (hV)\ave - k U\ave + F\ave,\label{eq:dUdt}
\end{equation}
where
\begin{equation*}(hV)\ave=\frac1{2\pi^2}\iint hV\dif x\dif y
\end{equation*} is the average topographic drag \cite{Davey1980}. If the average zonal velocity $U\ave$ is prescribed to be a constant, the average zonal force will be \begin{equation*}F\ave = k U\ave - (hV)\ave,
\end{equation*} which may vary with time. Therefore, a specified average zonal velocity gives the average zonal force implicitly, making the equations of the vertical vorticity well-posed.

The additional condition can even be given as
\begin{equation*}a F\ave + b U\ave + c=0,
\end{equation*}
 where $a$, $b$, and $c$ are constants, as in the study of the plane Poiseuille flow \cite{Barkley1990}. This additional condition includes the constant $F\ave$ condition and the constant $U\ave$ condition as special cases, and also permits a class of more general conditions
 \begin{equation*}F\ave = -\frac ba U\ave - \frac ca\end{equation*} when $a,b\neq 0$. In Section \ref{sec:discussion}, we will show that Charney and Devore solved the equations of the vertical vorticity (\ref{eq:psi_final})--(\ref{eq:psi_bc_final}) under this type of additional condition, although they did not realise it. The additional condition
 \begin{equation*}F\ave = -\frac ba U\ave - \frac ca\end{equation*} seems to lack a clear physical meaning, but it is useful in the calculation of multiple equilibrium states that are unstable in both the constant $F\ave$ calculation and the constant $U\ave$ calculation.

In addition, the ill-posedness of the equations of the vertical vorticity can also be revealed with the uniqueness of the solution of Poisson's equation. The equations of the vertical vorticity $\zeta$ depend on the stream function $\psi$, so we need to solve the Poisson's equation
\begin{equation*}\nabla^2 \psi = \zeta\end{equation*} in each time step, which has a unique solution $\psi_0$ under the boundary conditions
\begin{equation*}\psi|_{x=0} = \psi|_{x=2\pi}\ \mbox{ and } \ \psi|_{y=0}=\psi|_{y=\pi}=0.\end{equation*} Therefore, the general solution of the Poisson's equation  under the boundary condition (\ref{eq:psi_bc_final}) is $\psi=\psi_0 + C_0 + C_1 y $, where $C_0$ and $C_1$ are constants to be determined. The constant $C_0$ can be chosen arbitrarily, but $C_1 y$ contributes to the term $\p \psi/\p y$ in the equation (\ref{eq:psi_final}) and should be specified. Specifying $C_1$ is equivalent to specifying the average zonal velocity $U\ave$ because
\begin{equation*}C_1 = \frac1\pi (\psi|_{y=\pi} - \psi|_{y=0}) = -U\ave.\end{equation*} Therefore, an additional condition on the average zonal velocity guarantees a unique solution of Poisson's equation of the stream function, and makes the equations of the vertical vorticity well-posed.

In this paper, we concentrate on a simple case with
\begin{equation} F'=kC\left(\sin y - \frac2\pi\right),
\label{n1}
\end{equation} where $C$ is a constant as in \cite{Charney1979}. Cases with constant average zonal forces and cases with constant average zonal velocities are both considered in the following.

\section{Nonlinear stability}
\label{sec:theorem}
In this section, we will prove that there is an asymptotically stable solution to the equations (\ref{eq:continuity_final})--(\ref{eq:bc_final}) under the condition of constant average zonal force or constant average zonal velocity for $F'$ given by (\ref{n1}) and  $h = C\eta \cos y$, where $C$ and $\eta>-1$ are constants.

When the average zonal force $F\ave$ is a constant, and $F'$ and $h$ are functions of only the meridional coordinate $y$, the basic flow of the equations (\ref{eq:continuity_final})--(\ref{eq:bc_final}) is a parallel flow:
\begin{equation}
U_0=\frac{1}{k} (F\ave + F'),\qquad V_0=0, \qquad P_0 = -\frac{1}{k} \int_0^y ( h(y')+\beta y') (F\ave + F'(y')) \ \dif y' + P_{00},\label{eq:basic_flow}
\end{equation}
where $P_{00}$ is any constant.

From the equations (\ref{eq:continuity_final})--(\ref{eq:bc_final}), the governing equations and boundary conditions of the disturbance $(u,v,p)=(U,V,P)-(U_0,0,P_0)$ are
\begin{align}
& \frac{\p u}{\p x} + \frac{\p v}{\p y} = 0,\label{eq:continuity_dist}\\
& \frac{\p u}{\p t} + U_0 \frac{\p u}{\p x} + v \frac{\dif U_0}{\dif y} + u \frac{\p u}{\p x} + v \frac{\p u}{\p y} - (h + \beta y)v = -\frac{\p p}{\p x} -ku,\label{eq:momentum_x_dist}\\
& \frac{\p v}{\p t} + U_0 \frac{\p v}{\p x} + u \frac{\p v}{\p x} + v \frac{\p v}{\p y} + (h + \beta y)u = -\frac{\p p}{\p y} -kv,\label{eq:momentum_y_dist}\\
& (u,v,p)(x+2\pi,y,t) = (u,v,p)(x,y,t) , \qquad v|_{y=0} = v|_{y=\pi} = 0,\label{eq:bc_dist}
\end{align}
where $(u,v)$ are the velocity components of the disturbance, and $p$ is the pressure perturbation.

Introducing the disturbance vertical vorticity $\zeta' = \p v/\p x - \p u/\p y$, we have
\begin{equation*}
\frac{\p \zeta'}{\p t} + U_0 \frac{\p \zeta'}{\p x} + u \frac{\p \zeta'}{\p x} + v \frac{\p \zeta'}{\p y} + (\frac{\dif h}{\dif y} - \frac{\dif^2 U_0}{\dif y^2} + \beta) v = -k\zeta'.
\end{equation*}
Then we have
\begin{align}
& \frac{\dif}{\dif t}\int_0^\pi \int_0^{2\pi} \frac{u^2+v^2}{2} \ \dif x \dif y = -\int_0^\pi \int_0^{2\pi} uv\frac{\dif U_0}{\dif y} \ \dif x \dif y - k \int_0^\pi \int_0^{2\pi} (u^2+v^2) \ \dif x \dif y,\label{eq:dEdt}\\
& \frac{\dif}{\dif t}\int_0^\pi \int_0^{2\pi} \frac{\zeta'^2}{2} \ \dif x \dif y = -\int_0^\pi \int_0^{2\pi} uv (\frac{\dif^2 h}{\dif y^2} - \frac{\dif^3 U_0}{\dif y^3}) \ \dif x \dif y - k \int_0^\pi \int_0^{2\pi} \zeta'^2 \ \dif x \dif y,\label{eq:dzetadt}
\end{align}
where integration by parts and the boundary condition (\ref{eq:bc_dist}) are used.

If we assume $F'=kC(\sin y - 2/\pi)$ and $F\ave = k (C' + 2C/\pi)$, where $C$ and $C'$ are constants, then we have $U_0 = C \sin y + C'$ according to (\ref{eq:basic_flow}). When
$h = C\eta \cos y$, where $\eta$ is a constant, (\ref{eq:dEdt})--(\ref{eq:dzetadt}) lead to
\begin{equation}
\frac{\dif}{\dif t}\int_0^\pi \int_0^{2\pi} (\zeta'^2 + (\eta-1)(u^2+v^2)) \ \dif x \dif y =  - 2k \int_0^\pi \int_0^{2\pi} (\zeta'^2 + (\eta-1)(u^2+v^2)) \ \dif x \dif y.\label{eq:remarkable_identity}
\end{equation}
which is an extension of the remarkable identity in the case of flat topography \cite{Majda2006}.

When $\eta > 1$, (\ref{eq:remarkable_identity}) implies that the disturbance kinetic energy $0.5\iint (u^2+v^2)\ \dif x\dif y$ will exponentially decay. However, we are more interested in the stability of the basic flow when the topography is flat ($\eta=0$). It is natural to ask when is the integral $\int_0^\pi \int_0^{2\pi} (\zeta'^2 + (\eta-1)(u^2+v^2)) \ \dif x \dif y$ positive definite. Actually, we have the following lemma:
\begin{lem}\label{lem:zeta2_u2v2}
If $u,v\in H^1((0,2\pi)\times(0,\pi))$ are real functions with period $2\pi$ in the zonal direction, and satisfy
\begin{align}
& \frac{\p u}{\p x} + \frac{\p v}{\p y} = 0,\label{eq:condition_continuity}\\
& v|_{y=0} = v|_{y=\pi} = 0,\label{eq:condition_bc}\\
& \int_0^{2\pi} u(x,0) \ \dif x = \int_0^{2\pi} u(x,\pi) \ \dif x,\qquad \int_0^\pi \int_0^{2\pi} u \ \dif x \dif y = 0,\label{eq:condition_u}
\end{align}
then
\begin{equation}
\int_0^\pi \int_0^{2\pi} \zeta'^2 \ \dif x \dif y \geq 2 \int_0^\pi \int_0^{2\pi} (u^2+v^2) \ \dif x \dif y.\label{eq:zeta2_u2v2}
\end{equation}
\end{lem}

To prove Lemma \ref{lem:zeta2_u2v2}, we need the following two versions of Wirtinger's inequality \cite{Beesack1958}:
\begin{lem}\label{lem:Wirtinger1}
If $f \in H^1(0,L_y)$, $f(0)=f(L_y)$, and $\int_0^{L_y} f \ \dif y = 0$, then
\begin{equation*}
\int_0^{L_y} (\frac{\dif f}{\dif y})^2 \ \dif y \geq (\frac{2\pi}{L_y})^2 \int_0^{L_y} f^2\ \dif y.
\end{equation*}
\end{lem}

\begin{lem}\label{lem:Wirtinger2}
If $f \in H^2(0,L_y)$ and $f(0)=f(L_y)=0$, then
\begin{equation*}
\int_0^{L_y} (\frac{\dif^2 f}{\dif y^2})^2 \ \dif y \geq (\frac{\pi}{L_y})^2 \int_0^{L_y} f^2\ \dif y.
\end{equation*}
\end{lem}

\begin{proof}[\bf{Proof of Lemma \ref{lem:zeta2_u2v2}}]
The Fourier expansions of $u$, $v$, and $\zeta'$ are
\begin{equation*}
u(x,y) = \sum_{m=-\infty}^{+\infty} \hat{u}_m(y) \me^{\mi m x},\quad
v(x,y) = \sum_{m=-\infty}^{+\infty} \hat{v}_m(y) \me^{\mi m x},\quad
\zeta'(x,y) = \sum_{m=-\infty}^{+\infty} (\mi m \hat{v}_m - \frac{\dif \hat{u}_m}{\dif y}) \me^{\mi m x},
\end{equation*}
where $\hat{u}_{-m}=\bar{\hat{u}}_m$ and $\hat{v}_{-m}=\bar{\hat{v}}_m$, and the overlines denote complex conjugations. $\hat{u}_0,\hat{v}_0\in H^1(0,\pi)$ are real functions. Substituting the Fourier expansions into (\ref{eq:condition_continuity})--(\ref{eq:condition_u}), we have
\begin{align}
& \mi m \hat{u}_m + \frac{\dif \hat{v}_m}{\dif y} = 0,\quad (m\geq 0),\label{eq:Fourier_condition_continuity}\\
& \hat{v}_m (0) = \hat{v}_m (\pi) = 0,\quad (m\geq 0),\label{eq:Fourier_condition_bc}\\
& \hat{u}_0 (0) = \hat{u}_0 (\pi) ,\qquad \int_0^\pi \hat{u}_0 \ \dif y = 0.\label{eq:Fourier_condition_u}
\end{align}
Then we have $\hat{v}_0=0$, $\hat{u}_m \in H^1(0,\pi)$, and $\hat{v}_m \in H^2(0,\pi)$ for $m\geq 1$, and
\begin{align}
& \frac{1}{2\pi} \int_0^\pi \int_0^{2\pi} (u^2+v^2) \ \dif x \dif y = \int_0^\pi \hat{u}_0^2 \ \dif y + \sum_{m\geq 1} \int_0^\pi \Big( \frac{2}{m^2} \Big| \frac{\dif \hat{v}_m}{\dif y} \Big|^2 + 2 | \hat{v}_m |^2 \Big) \ \dif y, \label{eq:dEdt_Fourier}\\
& \frac{1}{2\pi} \int_0^\pi \int_0^{2\pi} \zeta'^2 \ \dif x \dif y = \int_0^\pi ( \frac{\dif \hat{u}_0}{\dif y} )^2 \ \dif y + \sum_{m\geq 1} \int_0^\pi \Big( \frac{2}{m^2} \Big| \frac{\dif^2 \hat{v}_m}{\dif y^2} \Big|^2 + 4 \Big| \frac{\dif \hat{v}_m}{\dif y} \Big|^2 + 2m^2 | \hat{v}_m |^2  \Big) \ \dif y.\label{eq:dzetadt_Fourier}
\end{align}
Using Lemma \ref{lem:Wirtinger1} for $\hat{u}_0$, we have
\begin{equation}
\int_0^\pi ( \frac{\dif \hat{u}_0}{\dif y} )^2 \ \dif y \geq 4 \int_0^\pi \hat{u}_0^2 \ \dif y.\label{eq:inequality_u0}
\end{equation}
Using Lemma \ref{lem:Wirtinger2} for $\hat{v}_m$, we have
\begin{align}
\int_0^\pi \Big( \frac{2}{m^2} \Big| \frac{\dif^2 \hat{v}_m}{\dif y^2} \Big|^2 + 4 \Big| \frac{\dif \hat{v}_m}{\dif y} \Big|^2 + 2m^2 | \hat{v}_m |^2 \Big) \ \dif y & \geq \int_0^\pi  \Big( 4 \Big| \frac{\dif \hat{v}_m}{\dif y} \Big|^2 + (\frac{2}{m^2}+2m^2) | \hat{v}_m |^2 \Big) \ \dif y \nonumber\\
& \geq 2\int_0^\pi  \Big( \frac{2}{m^2} \Big| \frac{\dif \hat{v}_m}{\dif y} \Big|^2 + 2 | \hat{v}_m |^2  \Big) \ \dif y\label{eq:inequality_vn}
\end{align}
for $m\geq 1$. Substituting (\ref{eq:inequality_u0}) and (\ref{eq:inequality_vn}) into (\ref{eq:dEdt_Fourier}) and (\ref{eq:dzetadt_Fourier}), we arrive at the inequality (\ref{eq:zeta2_u2v2}), which completes the proof.
\end{proof}

\begin{rmk}
The conclusion of Lemma \ref{lem:zeta2_u2v2} is optimal. When $(u,v) = (\sin x \cos y, -\cos x \sin y)$, $$\int_0^\pi \int_0^{2\pi} \zeta'^2 \ \dif x \dif y = 2 \int_0^\pi \int_0^{2\pi} (u^2+v^2) \ \dif x \dif y.$$
\end{rmk}

Under the conditions (\ref{eq:condition_continuity})--(\ref{eq:condition_u}), the integral $\int_0^\pi \int_0^{2\pi} (\zeta'^2 + (\eta-1)(u^2+v^2)) \ \dif x \dif y$ is positive definite for $\eta > -1$. A solution to the equations of the disturbance (\ref{eq:continuity_dist})--(\ref{eq:bc_dist}), however, does not have to satisfy the condition (\ref{eq:condition_u}) for all $t\geq 0$. According to the equation (\ref{eq:momentum_x_dist}), we have
\begin{align}
& \frac{\dif}{\dif t} \int_0^{2\pi} (u(x,0)-u(x,\pi)) \ \dif x = -k \int_0^{2\pi} (u(x,0)-u(x,\pi)) \ \dif x, \label{eq:decay_boundary}\\
& \frac{\dif}{\dif t} \int_0^\pi \int_0^{2\pi} u \ \dif x \dif y = -k \int_0^\pi \int_0^{2\pi} u \ \dif x \dif y.\label{eq:decay_u_ave}
\end{align}
Although the integrals $\int_0^{2\pi} (u(x,0)-u(x,\pi)) \ \dif x$ and $\int_0^\pi \int_0^{2\pi} u \ \dif x \dif y$ will not vanish if they do not vanish at $t=0$, they will both decay exponentially with time as $\exp(-kt)$. Their influence on the conclusion of Lemma \ref{lem:zeta2_u2v2} is characterised by the following lemma:

\begin{lem}\label{lem:zeta2_u2v2_u0}
If $u,v\in H^1((0,2\pi)\times(0,\pi))$ are real functions with period $2\pi$ in the zonal direction, and satisfy (\ref{eq:condition_continuity}) and (\ref{eq:condition_bc}), then
\begin{align}
\int_0^\pi \int_0^{2\pi} \zeta'^2 \ \dif x \dif y \geq & 2 \int_0^\pi \int_0^{2\pi} (u^2+v^2) \ \dif x \dif y - \frac{2}{\pi^2} \Big( \int_0^\pi \int_0^{2\pi} u \ \dif x \dif y \Big) ^2
\nonumber \\
&- \frac{\pi^2-3}{6\pi^2} \Big( \int_0^{2\pi} (u(x,0)-u(x,\pi)) \ \dif x \Big) ^2.\label{eq:zeta2_u2v2_u0}
\end{align}
\end{lem}

\begin{proof}[\bf{Proof}]
Performing the same Fourier expansions as in the proof of Lemma \ref{lem:zeta2_u2v2}, we still have
(\ref{eq:dEdt_Fourier}), (\ref{eq:dzetadt_Fourier}), and the inequality (\ref{eq:inequality_vn}). Using Lemma \ref{lem:Wirtinger1} for
\begin{equation*}
f(y)=\hat{u}_0(y) - \frac{1}{\pi}\int_0^{\pi} \hat{u}_0(y') \ \dif y' + (\frac{y}{\pi} - \frac{1}{2}) (\hat{u}_0(0) - \hat{u}_0(\pi)),
\end{equation*}
we have
\begin{align*}
& \int_0^\pi ( \frac{\dif \hat{u}_0}{\dif y} )^2 \ \dif y - \frac{1}{\pi} (\hat{u}_0(0) - \hat{u}_0(\pi))^2\\
& \geq 4 \int_0^\pi \hat{u}_0^2 \ \dif y + 8 (\hat{u}_0(0) - \hat{u}_0(\pi))  \int_0^{\pi} (\frac{y}{\pi} - \frac{1}{2})\hat{u}_0 \ \dif y  - \frac{4}{\pi} \Big( \int_0^{\pi} \hat{u}_0 \ \dif y \Big) ^2 + \frac{\pi}{3}(\hat{u}_0(0) - \hat{u}_0(\pi))^2\\
& \geq 2 \int_0^\pi \hat{u}_0^2 \ \dif y - \frac{4}{\pi} \Big( \int_0^{\pi} \hat{u}_0 \ \dif y \Big) ^2 - \frac{\pi}{3}(\hat{u}_0(0) - \hat{u}_0(\pi))^2,
\end{align*}
which is equivalent to
\begin{equation}
\int_0^\pi ( \frac{\dif \hat{u}_0}{\dif y} )^2 \ \dif y \geq 2 \int_0^\pi \hat{u}_0^2 \ \dif y - \frac{1}{\pi^3} \Big( \int_0^\pi \int_0^{2\pi} u \ \dif x \dif y \Big) ^2 - \frac{\pi^2-3}{12\pi^3} \Big( \int_0^{2\pi} (u(x,0)-u(x,\pi)) \ \dif x \Big) ^2.\label{eq:inequality_u0_new}
\end{equation}
Substituting (\ref{eq:inequality_u0_new}) and (\ref{eq:inequality_vn}) into (\ref{eq:dEdt_Fourier}) and (\ref{eq:dzetadt_Fourier}), we obtain the inequality (\ref{eq:zeta2_u2v2_u0}), which completes the proof.
\end{proof}

\begin{thm}\label{thm:Fave}
If $F' = k C (\sin y - 2/\pi)$, $F\ave = k (C' + 2C/\pi)$, and $h = C\eta \cos y$, where $C$, $C'$, and $\eta$ are constants, and $\eta > -1$, the basic flow $\bm{U}_0 = (C \sin y + C',0)$ is stable in the sense that
\begin{equation}
\| \bm{U}(t)-\bm{U}_0 \|_{H^1} \leq \Big( \frac{3}{1+\eta'} C_1 + \frac{2(2-\eta')}{\pi^2(1+\eta')} C_2 + \frac{(\pi^2-3)(2-\eta')}{6\pi^2(1+\eta')} C_3 \Big) ^{1/2} \exp(-kt), \label{eq:thm:Fave}
\end{equation}
where
\begin{align*}
& \eta' = \min\{\eta, 2\}, \\
& C_1 = \Big( \int_0^\pi \int_0^{2\pi} (\zeta'^2 + (\eta-1)(u^2+v^2)) \ \dif x \dif y \Big)_{t=0}, \\
& C_2 = \Big( \int_0^\pi \int_0^{2\pi} u \ \dif x \dif y \Big)^2_{t=0}, \\
& C_3 = \Big( \int_0^{2\pi} (u(x,0)-u(x,\pi)) \ \dif x \Big)^2_{t=0}.
\end{align*}
\end{thm}

\begin{proof}[\bf{Proof}] When $-1 < \eta < 2$, we have the following estimate of the $H^1$ norm of $\bm{u}=(u,v)$:
\begin{align*}
\| \bm{u} \|_{H^1}^2 = & \int_0^\pi \int_0^{2\pi} (\zeta'^2 + u^2+v^2) \ \dif x \dif y \\
\leq & \frac{3}{1+\eta} \int_0^\pi \int_0^{2\pi} (\zeta'^2 + (\eta-1)(u^2+v^2)) \ \dif x \dif y + \frac{2(2-\eta)}{\pi^2(1+\eta)} \Big( \int_0^\pi \int_0^{2\pi} u \ \dif x \dif y \Big) ^2\\
& + \frac{(\pi^2-3)(2-\eta)}{6\pi^2(1+\eta)} \Big( \int_0^{2\pi} (u(x,0)-u(x,\pi)) \ \dif x \Big)^2.
\end{align*}
When $\eta \geq 2$, it is obvious that
\begin{equation*}
\| \bm{u} \|_{H^1}^2 \leq \int_0^\pi \int_0^{2\pi} (\zeta'^2 + (\eta-1)(u^2+v^2)) \ \dif x \dif y.
\end{equation*}
Furthermore, noticing (\ref{eq:remarkable_identity}), (\ref{eq:decay_boundary}), and (\ref{eq:decay_u_ave}), we  finish the proof of  the  theorem.
\end{proof}

In \cite{Chen2016b}, we proved a weaker version of Theorem \ref{thm:Fave}, which stated that $\| \zeta' \|_{L^2}$ would decay exponentially with time provided that $\eta \geq 0$. But the proof therein is not strict. The governing equations in \cite{Chen2016b} need an additional condition on the average zonal force or the average zonal velocity to be well-posed, as in the work of Charney and DeVore \cite{Charney1979}.

Now we consider the case with a constant average zonal velocity instead of a constant average zonal force. When the average zonal velocity $U\ave$ is a constant, and $h$ is a function of only the meridional coordinate $y$, the average zonal force $F\ave = k U\ave - (hV)\ave = k U\ave$ is also a constant according to the equation (\ref{eq:dUdt}), where $(hV)\ave=(\p (h\psi)/\p x )\ave = 0$ is used. Therefore we have
\begin{thm}\label{thm:Uave}
If $F' = k C (\sin y - 2/\pi)$, $U\ave = C' + 2C/\pi$, and $h = C\eta \cos y$, where $C$, $C'$, and $\eta$ are constants, and $\eta > -1$, the basic flow $\bm{U}_0 = (C \sin y + C',0)$ is stable in the sense that the inequality (\ref{eq:thm:Fave}) holds.
\end{thm}

\section{Numerical methods}
\label{sec:method}
\subsection{Direct numerical simulation}
\label{sec:DNS}
The periodic boundary conditions for $U$, $V$, and $P$ in the zonal direction permit Fourier expansions of these primitive variables. Denoting the Fourier expansion in the zonal direction of a real function $f(x,y,t)$ as $\sum_m  \hat{f}_m(y,t) \me ^{\mi mx}$ or $\sum_m  (\mathcal{F}f)_m \me ^{\mi mx}$, where $\hat{f}_0 \in \mathbb{R}$ and $\hat{f}_{-m}=\bar{\hat{f}}_{m}$, with $\bar{\hat{f}}_{m}$ representing the complex conjugation of $\hat{f}_{m}$. The equations of primitive variables (\ref{eq:continuity_final})--(\ref{eq:bc_final}) with constant average zonal force $F\ave$ lead to
\begin{align}
& \mi m \hat{U}_m + \frac{\p \hat{V}_m}{\p y} = 0,\\
& \frac{\p \hat{U}_m}{\p t} + \mathcal{F}(U \frac{\p U}{\p x} + V \frac{\p U}{\p y})_m - \mathcal{F}(hV)_m - \beta y\hat{V}_m = -\mi m \hat{P}_m -k\hat{U}_m + \delta_{m,0} F\ave + \hat{F}'_m,\\
& \frac{\p \hat{V}_m}{\p t} + \mathcal{F}(U \frac{\p V}{\p x} + V \frac{\p V}{\p y})_m + \mathcal{F}(hU)_m + \beta y\hat{U}_m = -\frac{\p \hat{P}_m}{\p y} -k\hat{V}_m,\\
& \hat{V}_m|_{y=0} = \hat{V}_m|_{y=\pi} = 0,
\end{align}
where $\delta_{m,0}=1$ for $m=0$, and $\delta_{m,0}=0$ for $1 \leq m\leq M-1$. $M$ is the number of independent Fourier modes used in the discretisation. Eliminating $\hat{P}_m$ from the above equations, we have
\begin{align}
& \frac{\p \hat{U}_0}{\p t} + \mathcal{F}(U \frac{\p U}{\p x} + V \frac{\p U}{\p y})_0 - \mathcal{F}(hV)_0 = -k\hat{U}_0 + F\ave  + \hat{F}'_0,\label{eq:U0hat}\\
& \hat{V}_0 = 0,
\end{align}
and
\begin{align}
& \hat{U}_m = \frac{\mi}{m} \frac{\p \hat{V}_m}{\p y},\\
 &\frac{\p}{\p t} (\frac{\p^2}{\p y^2} - m^2) \hat{V}_m + \mi m \mathcal{F}(U \nabla^2 V - V \nabla^2 U)_m + \mi m \mathcal{F}(U \frac{\p h}{\p x} + V \frac{\p h}{\p y})_m + \mi m \beta \hat{V}_m
 \nonumber\\ &=-k (\frac{\p^2}{\p y^2} - m^2) \hat{V}_m - \mi m \frac{\p \hat{F}'_m}{\p y} ,\\
& \hat{V}_m|_{y=0} = \hat{V}_m|_{y=\pi} = 0,\label{eq:Vmhat_bc}
\end{align}
for $1 \leq m \leq M-1$.

If the average zonal velocity $U\ave$ is prescribed as a constant, the average zonal force will be $F\ave = k U\ave- (hV)\ave$ according to the equation (\ref{eq:dUdt}). Instead of the equation (\ref{eq:U0hat}), we have
\begin{equation}
\frac{\p \hat{U}_0}{\p t} + \mathcal{F}(U \frac{\p U}{\p x} + V \frac{\p U}{\p y})_0 - \mathcal{F}(hV)_0 = -k\hat{U}_0 + k U\ave - (hV)\ave  + \hat{F}'_0.
\end{equation}

If the additional condition is $a F\ave + b U\ave + c = 0$, where $a$, $b$, and $c$ are constants, and $a,b\neq 0$, we have the following equation instead of the equation (\ref{eq:U0hat}):
\begin{equation}
\frac{\p \hat{U}_0}{\p t} + \mathcal{F}(U \frac{\p U}{\p x} + V \frac{\p U}{\p y})_0 - \mathcal{F}(hV)_0 = -k\hat{U}_0 - \frac{b}{a} U\ave - \frac{c}{a}   + \hat{F}'_0.
\end{equation}

The Chebyshev-Tau method with $N$ Chebyshev polynomials is used for the discretisation in the meridional direction. In the temporal discretisation, linear terms except terms containing $h$ are discretised implicitly with the second-order Crank-Nicolson scheme, and other terms are discretised explicitly with the four-stage third-order Runge-Kutta scheme \cite{SIMSON}.

\subsection{Pseudo-arclength continuation method}
\label{sec:continuation}
Direct numerical simulations are difficult to obtain linearly unstable equilibrium states, because any small-amplitude disturbances will lead the solutions away from the unstable equilibrium states in the phase space. Therefore, it is more adequate to solve the nonlinear steady equations directly, instead of calculating equilibrium states as final states after a sufficiently long time in direct numerical simulations. The nonlinear steady equations are solved with a pseudo-arclength continuation method.

Removing the time derivative terms in (\ref{eq:U0hat})--(\ref{eq:Vmhat_bc}), expanding the nonlinear terms, and eliminating $\hat{U}_m$, we have
\begin{equation}
-\sum_{m\geq 1} ( \frac{\mi}{m} \bar{\hat{V}}_m \Dif^2 \hat{V}_m - \frac{\mi}{m}  \hat{V}_m \Dif^2 \bar{\hat{V}}_m ) + \sum_{m\geq 1} ( \bar{\hat{h}}_m \hat{V}_m + \hat{h}_m \bar{\hat{V}}_m ) - k\hat{U}_0 + F\ave  + \hat{F}'_0 = 0, \label{eq:continuation_U0}
\end{equation}
and
\begin{align}
0=&\sum_{\substack{m_1+m_2=m \\ m_1\neq 0}} \Bigg( -\frac{m}{m_1} (\Dif \hat{V}_{m_1}) (\Dif^2-m_2^2) \hat{V}_{m_2} + \frac{m}{m_1} \hat{V}_{m_2} (\Dif^2-m_1^2) \Dif \hat{V}_{m_1}
\nonumber\\
&\qquad\qquad\qquad-\frac{\mi m m_2}{m_1} (\Dif \hat{V}_{m_1}) \hat{h}_{m_2} + \mi m \hat{V}_{m_1} \Dif \hat{h}_{m_2} \Bigg) \nonumber \\
& + (\mi m \hat{U}_0 + k) (\Dif^2-m^2) \hat{V}_m - \mi m (\Dif^2 \hat{U}_0 - \beta) \hat{V}_m  - m^2 \hat{U}_0 \hat{h}_m + \mi m \Dif \hat{F}'_m  \label{eq:continuation_V}
\end{align}
for $1 \leq m \leq M-1$, where $\Dif \equiv \dif /\dif y$.

Using the Chebyshev-Tau method, we expand $\hat{U}_0$ and $\hat{V}_m$ with $N$ Chebyshev polynomials, and denote the column vectors composed of the coefficients as $\hat{\bm{U}}_0$ and $\hat{\bm{V}}_m$, respectively. Then the equations (\ref{eq:continuation_U0}) and (\ref{eq:continuation_V}) with the boundary condition (\ref{eq:Vmhat_bc}) constitute a nonlinear equation $\bm{G}(\bm{f}, F\ave) = \bm{0}$, where
\begin{align*}
& \bm{f} = (\hat{\bm{U}}_0\tr \quad \R \hat{\bm{V}}_1\tr \quad \I \hat{\bm{V}}_1\tr \quad \cdots \quad \R \hat{\bm{V}}_{M-1}\tr \quad \I \hat{\bm{V}}_{M-1}\tr)\tr.
\end{align*}
The superscript ``$\mathrm{T}$'' denotes the transpose of a vector or a matrix. ``$\mathrm{Re}$'' and ``$\mathrm{Im}$'' represent the real part and the imaginary part of a complex vector, respectively.

We use the following predictor-corrector continuation method \cite{Keller1977, Keller1987, Seydel1988} to find a branch of solutions of $\bm{G}(\bm{f}, F\ave) = \bm{0}$.

\subsubsection{Predictor (the first solution)}
\label{sec:predictor1}
The initial guess of the first solution is chosen as an approximation solution of the equations (\ref{eq:continuation_U0}), (\ref{eq:continuation_V}), and (\ref{eq:Vmhat_bc}) at a sufficiently large $F\ave$. For example, if $h = h_0 \cos 2x \sin y$, then an approximation solution is $F\ave=0.003$, $\hat{U}_0 = (F\ave + \hat{F}'_0)/k$, $\hat{V}_2 = 0.2\mi  h_0 \sin y$, and $\hat{V}_m=0$ for $1\leq m \leq M-1$ and $m\neq 2$. Denote this approximation solution as $(\bm{f}^0_1, F_{\mathrm{ave},1})$ or $(\bm{f}^0_1, F_1)$. The subscript ``$\mathrm{ave}$'' will be omitted in the following for conciseness.

\subsubsection{Corrector (the first solution)}
The iterative scheme for the calculation of the first solution is
\begin{align*}
& \frac{\p \bm{G}}{\p \bm{f}}(\bm{f}^{n-1}_1, F_1) \Delta \bm{f}^{n}_1 = -\bm{G}(\bm{f}^{n-1}_1, F_1), \\
& \bm{f}^{n}_1 = \bm{f}^{n-1}_1 + \Delta \bm{f}^{n}_1,
\end{align*}
for $n\geq 1$.
If the sequence $\{\bm{f}^{n}_1\}_{n\geq 0}$ converges to $\bm{f}_1$, then the first solution is defined as $(\bm{f}_1, F_1)$. Otherwise, a more accurate approximation solution is needed in Section \ref{sec:predictor1}. For example, we can choose the same form of the approximation solution, but at a larger $F\ave$.

\subsubsection{Predictor (the $\nu$-th solution, $\nu \geq 2$)}
After the solution $(\bm{f}_{\nu-1}, F_{\nu-1})$ is obtained, we search the next solution $(\bm{f}_{\nu}, F_{\nu})$ which satisfies $||(\bm{f}_{\nu} - \bm{f}_{\nu-1}, F_{\nu} - F_{\nu-1})|| = \Delta s$, where $||(\bm{f}, F)||^2 = \bm{f}\tr \bm{W} \bm{f} + F^2/(2k^2)$. $\bm{W}$ is a weight matrix so that $\bm{f}\tr \bm{W} \bm{f}$ is the discretisation of the average kinetic energy
\begin{equation}
E\ave = \frac{1}{2\pi^2} \int_0^\pi \int_0^{2\pi} \frac{U^2+V^2}{2} \ \dif x \dif y. \label{eq:Ek}
\end{equation}

An approximation of the $\nu$-th solution can be obtained by a tangent predictor \cite{Keller1987}:
\begin{align*}
& \frac{\p \bm{G}}{\p \bm{f}}(\bm{f}_{\nu-1}, F_{\nu-1}) \Delta \tilde{\bm{f}}^{0}_{\nu} = -\frac{\p \bm{G}}{\p F}(\bm{f}_{\nu-1}, F_{\nu-1}), \\
& (\Delta \bm{f}^{0}_{\nu}, \Delta F^{0}_{\nu}) = \pm \frac{(\Delta \tilde{\bm{f}}^{0}_{\nu}, 1)}{||(\Delta \tilde{\bm{f}}^{0}_{\nu}, 1)||} \Delta s, \\
& (\bm{f}_{\nu}, F_{\nu}) = (\bm{f}_{\nu-1}, F_{\nu-1}) + (\Delta \bm{f}^{0}_{\nu}, \Delta F^{0}_{\nu}).
\end{align*}
The plus or minus sign is determined by demanding the angle between $(\Delta \bm{f}^{0}_{\nu}, \Delta F^{0}_{\nu})$ and $(\bm{f}_{\nu-1} - \bm{f}_{\nu-2}, F_{\nu-1} - F_{\nu-2})$ is less than $90^\circ$ when $\nu\geq 3$, where the angle is calculated with the innerproduct corresponding to the norm $||\cdot||$. When $\nu=2$, we choose the $(\Delta \bm{f}^{0}_2, \Delta F^{0}_2)$ with a negative $\Delta F^{0}_2$ to calculate the solutions $(f,F)$ with $F < F_1$. The other half of the branch of solutions with $F > F_1$ can also be obtained by choosing a positive $\Delta F^{0}_2$ in a new calculation.\

\subsubsection{Corrector (the $\nu$-th solution, $\nu \geq 2$)}
\label{sec:Corrector_nu}
The iterative scheme for the calculation of the $\nu$-th solution $(\bm{f}_{\nu}, F_{\nu})$ ($\nu \geq 2$) is
\begin{align*}
& \left(
\begin{matrix}
\frac{\p \bm{G}}{\p \bm{f}}(\bm{f}^n_{\nu}, F^n_{\nu}) & \frac{\p \bm{G}}{\p F}(\bm{f}^n_{\nu}, F^n_{\nu}) \\
2(\bm{f}^n_{\nu} - \bm{f}_{\nu-1})\tr \bm{W} & (F^n_{\nu} - F_{\nu-1})/k^2
\end{matrix}
\right)
\left(
\begin{matrix}
\Delta \bm{f}^{n+1}_{\nu} \\
\Delta F^{n+1}_{\nu}
\end{matrix}
\right)
=
\left(
\begin{matrix}
-\bm{G}(\bm{f}^n_{\nu}, F^n_{\nu}) \\
- A_{\nu}(\bm{f}^n_{\nu}, F^n_{\nu})
\end{matrix}
\right), \\
& (\bm{f}^{n+1}_{\nu}, F^{n+1}_{\nu}) = (\bm{f}^{n}_{\nu}, F^{n}_{\nu}) + (\Delta \bm{f}^{n+1}_{\nu}, \Delta F^{n+1}_{\nu}),
\end{align*}
for $n\geq 1$, where
\begin{equation*}
A_{\nu}(\bm{f}, F) = ||(\bm{f} - \bm{f}_{\nu-1}, F - F_{\nu-1})||^2 - \Delta s^2.
\end{equation*}
This scheme is inspired by the Taylor expansions of $\bm{G}(\bm{f}^{n+1}_{\nu}, F^{n+1}_{\nu}) = \bm{0}$ and $A_{\nu}(\bm{f}^{n+1}_{\nu}, F^{n+1}_{\nu}) = 0$ at $(\bm{f}^{n}_{\nu}, F^{n}_{\nu})$.

If the sequence $\{(\bm{f}^{n}_{\nu}, F^{n}_{\nu})\}_{n\geq 0}$ converges, the $\nu$-th solution $(\bm{f}_{\nu}, F_{\nu})$ is defined as the limit of the sequence. Otherwise, we change the step from $\Delta s = \Delta s_0$ to $\Delta s = 0.5\Delta s_0$, and recalculate the $\nu$-th solution.

The branch of solution ends where new solutions of $\bm{G}(\bm{f}, F\ave) = \bm{0}$ cannot be obtained no matter how small the step is, which means either there are no more equilibrium states or the equilibrium states need a better resolution (larger $M$ and $N$) to be calculated.

\section{Results}
\label{sec:results}
In Section \ref{sec:theorem}, we have proved the asymptotic stability of the basic flow for a zonal-invariant topography $h = C\eta \cos y$ in a region of parameter space. In the following, we perform numerical simulations for a more realistic topography $h = h_0 \cos 2x \sin y$, where $h_0$ is a constant. This topography was also considered by Charney and Devore \cite{Charney1979}. When $h = h_0 \cos 2x \sin y$ and $F'=\sqrt{2}k\psiAO (\sin y - 2/\pi)$, the equations (\ref{eq:continuity_final})--(\ref{eq:bc_final}) can be solved in the following three spaces:
\begin{align*}
&  \mathbb{S} = \{ (U,V) \in H^1 \times H^1 \ | \  \frac{\p U}{\p x} + \frac{\p V}{\p y} = 0,\quad V(x,0)=V(x,\pi)=0 \}, \\
& \mathbb{S}_{even} = \{(U,V) \in \mathbb{S} \ | \  U(x+\pi,y) = U(x,y), \quad V(x+\pi,y) = V(x,y) \}, \\
& \mathbb{S}_{sym} = \{(U,V) \in \mathbb{S}_{even}\ | \ U(x+\frac{\pi}{2}, \pi-y) = U(x,y), \quad V(x+\frac{\pi}{2}, \pi-y) = -V(x,y) \},
\end{align*}
where $H^1$ represents the space of real functions in $H^1((0,2\pi)\times(0,\pi))$ with period $2\pi$ in the zonal direction.

The equilibrium states of the equations (\ref{eq:continuity_final})--(\ref{eq:bc_final}) are calculated with the pseudo-arclength continuation method in the space $\mathbb{S}_{sym}$ for $k = 0.01$, $\beta = 0.25$, and $\psiAO = 0.2$. We use 128 even or odd Chebyshev polynomials in the meridional direction, and use 16 even Fourier modes in the zonal direction. It is equivalent to setting $N=256$ and $M=32$ in Section \ref{sec:continuation} because of the symmetry of the space $\mathbb{S}_{sym}$. The maximum step in the continuation method is 0.01.

\subsection{Multiple equilibrium states for $h_0=0.2$}
\begin{figure}[htbp]
\centering \includegraphics[width=80mm]{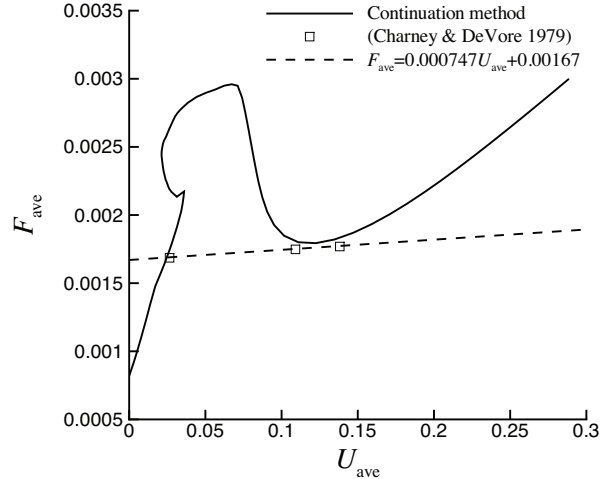}
\caption{The average zonal force and the average zonal velocity of the equilibrium states for $h_0 = 0.2$. The solid line is composed of the equilibrium states calculated with the continuation method. The three equilibrium states calculated by Charney and DeVore \cite{Charney1979} are shown by squares (Table \ref{tab:psi}), and satisfy $F\ave=0.000747U\ave+0.00167$ (Eqn (\ref{eq:additional_condition}) when $k=0.01$ and $\psiAO=0.2$), which is represented by a dashed line.}
\label{fig:Fave_Uave_h0}
\end{figure}

\begin{figure}[htbp]
\centering \includegraphics[width=158mm]{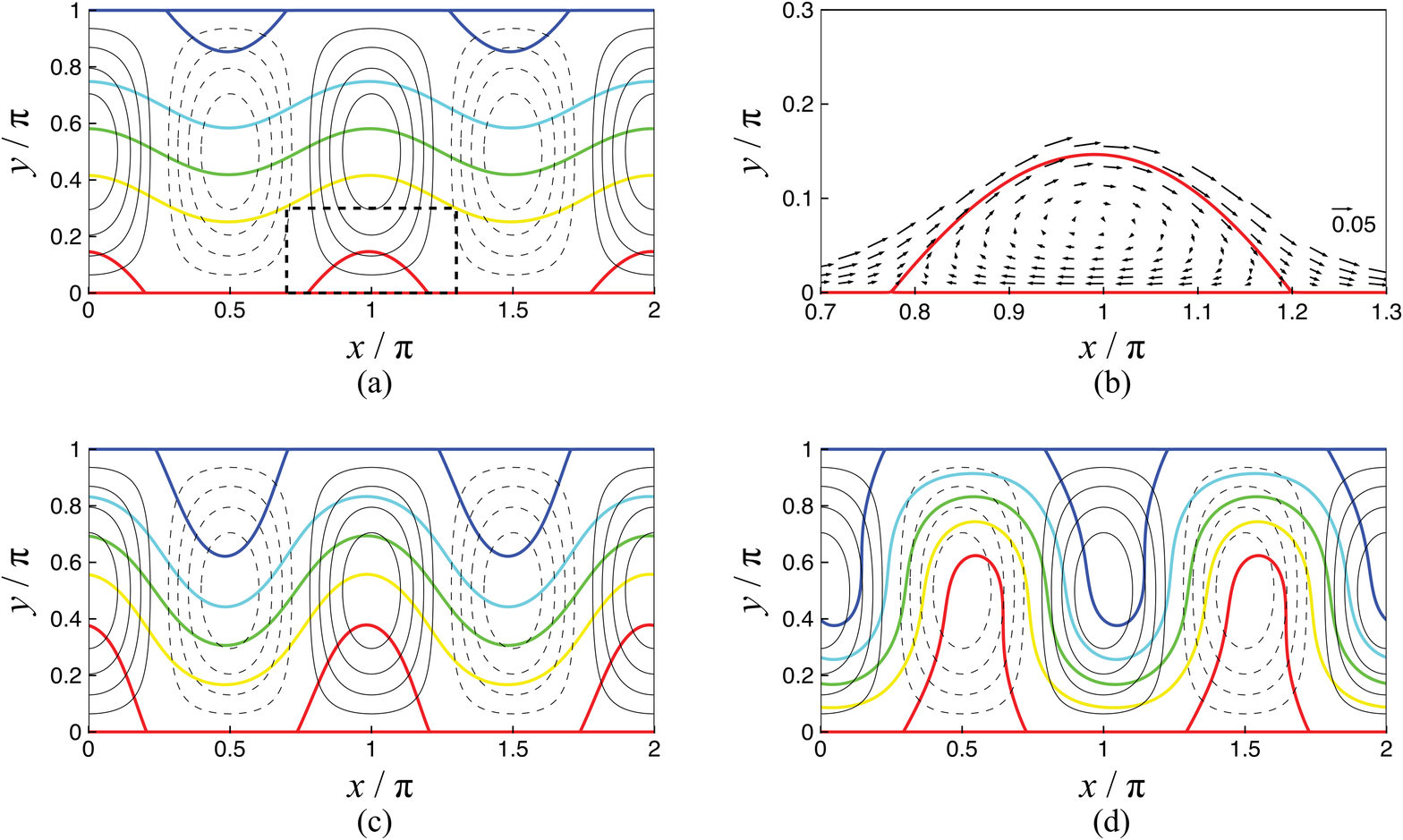}
\caption{The contours of the stream function $\psi$ (thick coloured lines) and the topography $h$ (thin black lines with dashed lines representing negative values) of the three equilibrium states for $h_0 = 0.2$ and $F\ave=0.002$. The average zonal velocity of these equilibrium states are $U\ave=0.170$ (a), $U\ave=0.0915$ (c), and $U\ave=0.0336$ (d). The thick coloured contours in each subfigure correspond to $\psi = 0$, $-0.25\pi U\ave$, $-0.5\pi U\ave$, $-0.75\pi U\ave$, and $-\pi U\ave$ from bottom to top. (b) is a close-up view of the dashed rectangle in (a) with velocity vectors plotted on a coarse mesh in the region $\psi > -0.005$, which shows a small-scale circulation attached to the boundary $y=0$.}
\label{fig:psi_F_0.002}
\end{figure}

When $h_0=0.2$, there are multiple equilibrium states with different $U\ave$ corresponding to the same $F\ave$ (Figure \ref{fig:Fave_Uave_h0}). For example, when $F\ave=0.002$, there are a ``high-index'' flow, a ``medium-index'' flow, and a ``low-index'' flow with $U\ave=0.170$, $0.0915$, and $0.0336$, respectively (Figure \ref{fig:psi_F_0.002}). The thick coloured contours in each subfigure of Figure \ref{fig:psi_F_0.002} correspond to $\psi = 0$, $-0.25\pi U\ave$, $-0.5\pi U\ave$, $-0.75\pi U\ave$, and $-\pi U\ave$ from bottom to top. In this paper, we always choose $\psi |_{y=0}=0$, and then we have $\psi |_{y=\pi}=-\pi U\ave$. The thin black lines are the contours of the topography $h = h_0 \cos 2x \sin y$ for reference, with dashed lines representing negative values.

There is a small-scale circulation in the region circled by the closed streamline $\psi=0$ in the dashed rectangle of Figure \ref{fig:psi_F_0.002}(a), which can be seen clearly in Figure \ref{fig:psi_F_0.002}(b). The stream function achieves its maximum $\psi_{\max}$ at the stagnation point $(x,y)=(0.99\pi, 0.08\pi)$ in Figure \ref{fig:psi_F_0.002}(b). Therefore, the circumferential flux of the small-scale circulation in the region $\psi>0$ is $\phi\cir = \psi_{\max} - 0$, where a positive flux represents a clockwise circulation.

\begin{figure}[htbp]
\centering \includegraphics[width=158mm]{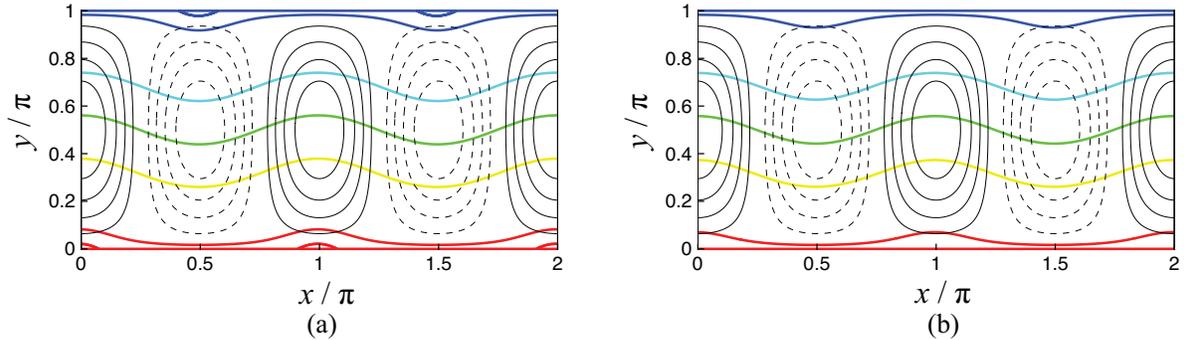}
\caption{The contours of the stream function $\psi$ (thick coloured lines) and the topography $h$ (thin black lines with dashed lines representing negative values) for $(F\ave, U\ave)=(0.00238, 0.220)$ (a) and $(F\ave, U\ave)=(0.00247, 0.230)$ (b) when $h_0 = 0.2$. The thick coloured contours in each subfigure correspond to $\psi = 0$, $-0.01\pi U\ave$, $-0.25\pi U\ave$, $-0.5\pi U\ave$, $-0.75\pi U\ave$, $-0.99\pi U\ave$, and $-\pi U\ave$ from bottom to top.}
\label{fig:psi_U_0.22_0.23}
\end{figure}

Actually, four regions contain small-scale circulations in Figure \ref{fig:psi_F_0.002}(a). Two of them attaching to the boundary $y=0$ have the same circumferential flux $\phi\cir$, and other two regions attaching to the boundary $y=\pi$ have the same circumferential flux $-\phi\cir$, because of the symmetry of the velocity field in the space $\mathbb{S}_{sym}$. The regions of small-scale circulations will shrink when $U\ave$ is increasing. When $U\ave \geq 0.23$, small-scale circulations no longer exist. The contours of the stream function when $U\ave=0.22$ and $0.23$ are plotted in Figure \ref{fig:psi_U_0.22_0.23}. In Figure \ref{fig:psi_U_0.22_0.23}(b), the contours $\psi=0$ and $\psi=-\pi U\ave$ coincide with the boundary $y=0$ and $y=\pi$, respectively.

In Charney and DeVore's numerical simulation, the low-index and the high-index equilibrium states are stable, whereas the medium-index equilibrium state is unstable \cite{Charney1979}. We examine the stabilities of some equilibrium states with direct numerical simulations. The velocity fields of the equilibrium states calculated with the pseudo-arclength continuation method are used as the initial velocity fields in direct numerical simulations. We use 192 even or odd Chebyshev polynomials and 192 even Fourier modes in the meridional direction and the zonal direction, respectively. It is equivalent to setting $N=M=384$ in Section \ref{sec:DNS} because of the symmetry of the space $\mathbb{S}_{sym}$. The time step is $\Delta t=0.02$.

\begin{figure}[htbp]
\centering \includegraphics[width=80mm]{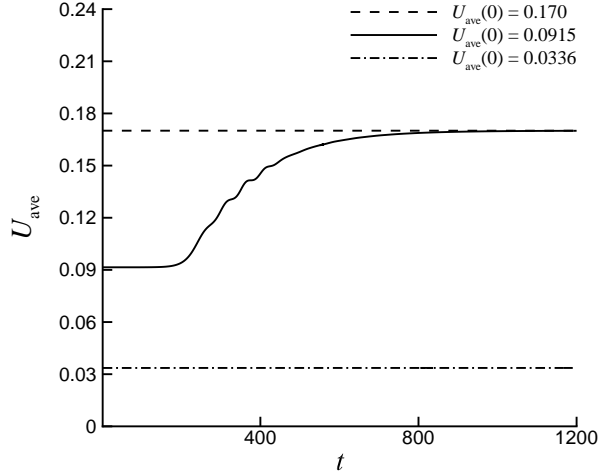}
\caption{The evolution of $U\ave$ from the three equilibrium states with $(F\ave, U\ave)=(0.002, 0.170)$ (dashed line), $(F\ave, U\ave)=(0.002, 0.0915)$ (solid line), and $(F\ave, U\ave)=(0.002, 0.0336)$ (dashed-dot line) in direct numerical simulations when $F\ave$ is fixed to be 0.002 and $h_0=0.2$.}
\label{fig:energy_Fave=0.002}
\end{figure}

When $F\ave=0.002$, the medium-index equilibrium state $(F\ave, U\ave)=(0.002, 0.0915)$ is unstable in $\mathbb{S}_{sym}$, and will evolve into the high-index equilibrium state $(F\ave, U\ave)=(0.002, 0.170)$ (Figure \ref{fig:energy_Fave=0.002}). The low-index and the high-index equilibrium states are stable in the direct numerical simulations.

\begin{figure}[htbp]
\centering {\includegraphics[width=158mm]{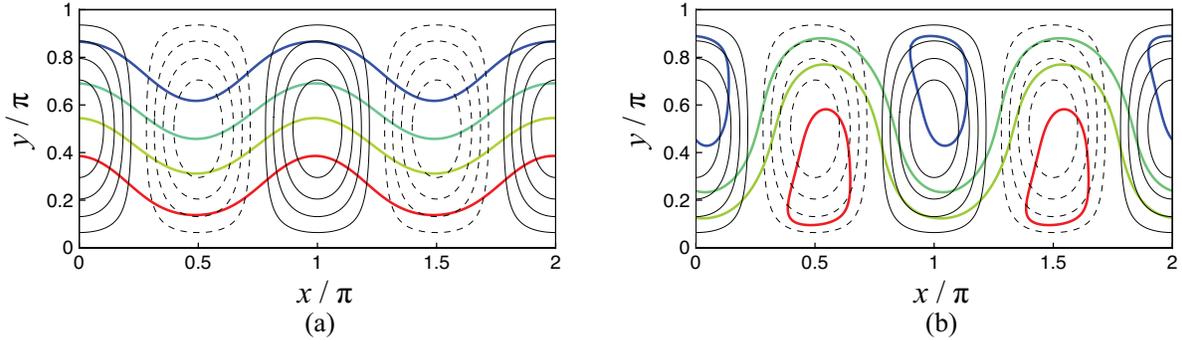}}
\caption{The contours of the stream function $\psi$ (thick coloured lines) and the topography $h$ (thin black lines with dashed lines representing negative values) of the ``high-index'' equilibrium state $U\ave=0.127$ (a) and the ``low-index'' equilibrium state $U\ave=0.0283$ (b) for $F\ave=0.0018$ and $h_0=0.2$, using the same contour levels as Fig.\,6 in \cite{Chen2016b}.}
\label{fig:psi_Fave=0.0018}
\end{figure}

In our previous numerical simulation \cite{Chen2016b}, we set the zonal component of the body force to be $F=\sqrt{2}k\psiAO\sin y$ and required the ``net pressure gradient'' to be zero, which was equivalent to setting constant average zonal force to be $F\ave=2\sqrt{2}k\psiAO/\pi$ and $F'=\sqrt{2}k\psiAO(\sin y - 2/\pi)$. Using the same parameters $k = 0.01$, $\beta = 0.25$, $h_0 = 0.2$, and $\psiAO = 0.2$, we obtain three equilibrium states for $F\ave=2\sqrt{2}k\psiAO/\pi \approx 0.0018$ with the present continuation method, including the two stable equilibrium states with $U\ave=0.127$ and $0.0283$ shown in Fig.\,6 of \cite{Chen2016b} and another equilibrium state with $U\ave=0.111$. We plot the contours of $\psi=-0.05$, $-0.15$, $-0.25$ and $-0.35$ for the equilibrium state $U\ave=0.127$ in Figure \ref{fig:psi_Fave=0.0018}(a), and plot the contours of $\psi=0.01$, $-0.03$, $-0.06$ and $-0.1$ for the equilibrium state $U\ave=0.0283$ in Figure \ref{fig:psi_Fave=0.0018}(b), which agree well with Fig.\,6 in \cite{Chen2016b}.

\subsection{Equilibrium states for $0\leq h_0 \leq 0.2$}
\begin{figure}[htbp]
\centering \includegraphics[width=100mm]{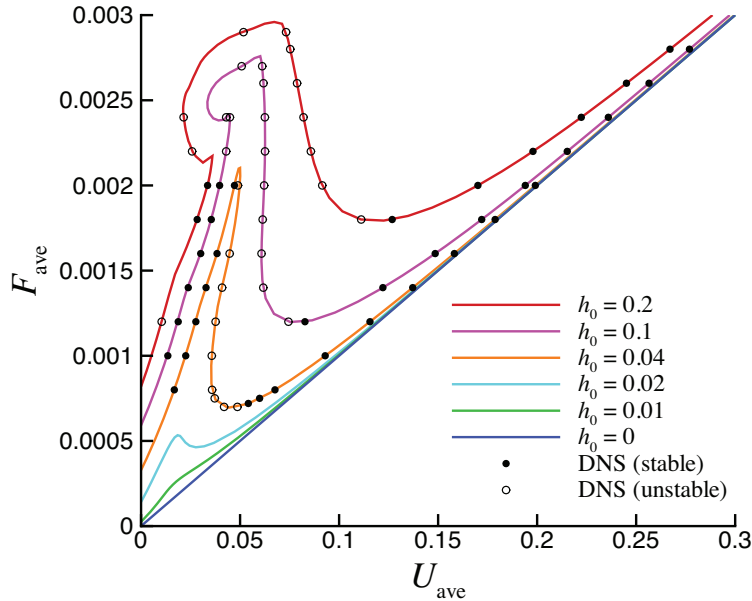}
\caption{The average zonal force and the average zonal velocity of the equilibrium states for $0 \leq h_0 \leq 0.2$. The solid lines are composed of the equilibrium states calculated with the continuation method. Some of these equilibrium states are verified with direct numerical simulations, and are plotted with solid circles and open circles, representing stable and unstable equilibrium states, respectively.}
\label{fig:Fave_Uave}
\end{figure}

When the topography is flat ($h_0=0$), the equilibrium state $U=(F\ave + F')/k$, $V=0$ is globally stable for each $F\ave$ according to Theorem \ref{thm:Fave}. The average zonal velocity is $U\ave = F\ave/k$ from the equation (\ref{eq:dUdt}). Actually, the same equilibrium state $U=U\ave + F'/k$, $V=0$ is also globally stable for each $U\ave$ according to Theorem \ref{thm:Uave}. So there is only one stable equilibrium state for each $F\ave$ or for each $U\ave$ when $h_0=0$. When $h_0=0.01$, there is only one equilibrium state for each $F\ave$ as in the case of $h_0=0$ (Figure \ref{fig:Fave_Uave}). When $h_0$ is increased to 0.02, there are three equilibrium states for each $F\ave$ when $0.00046 \leq F\ave \leq 0.00053$.

The stabilities of some equilibrium states are examined with direct numerical simulations. In Figure \ref{fig:Fave_Uave}, the stable and unstable equilibrium states are represented by solid and open circles, respectively. When $h_0 = 0.04$, the stabilities of the equilibrium states change when $F\ave$ is near its extreme values. When $h_0 = 0.1$ or $h_0 = 0.2$, however, the stabilities of the equilibrium states may be influenced by limit cycles which are not considered in the present continuation method. For example, when $h_0 = 0.1$, the equilibrium state with $(F\ave, U\ave)=(0.0022, 0.0429)$ will evolve into a limit cycle in the direct numerical simulation, with the average zonal velocity $U\ave$ oscillating between 0.0407 and 0.0448 (Figure \ref{fig:limit_cycle_Fave=0.0022}(a)). The projection of the phase space trajectory on the average kinetic energy--average zonal velocity ($E\ave$-$U\ave$) plane in Figure \ref{fig:limit_cycle_Fave=0.0022}(b) shows a stable limit cycle. The period of the limit cycle is $T\approx 56$ (Figure \ref{fig:limit_cycle_Fave=0.0022}(c)). The contours of the stream functions when $U\ave$ achieves its minimum and maximum are shown in Figure \ref{fig:limit_cycle_Fave=0.0022}(d) and (e), respectively. The small-scale circulations almost do not move in the zonal direction as a travelling wave, but they intrude further into the zonal flow in the meridional direction when $U\ave$ is minimum than when $U\ave$ is maximum.

\begin{figure}[htbp]
\centering \includegraphics[width=158mm]{Figure_7.eps}
\caption{The evolution of the equilibrium state with $(F\ave, U\ave)=(0.0022, 0.0429)$ in the direct numerical simulation when $F\ave$ is fixed to be 0.0022 and $h_0=0.1$. (a) Evolution of $U\ave$. (b) Phase space trajectory (projection on the $E\ave$-$U\ave$ plane). (c) Evolution of $U\ave$ and $E\ave$ (a close-up view). The contours of the stream function $\psi$ (thick coloured lines) and the topography $h$ (thin black lines with dashed lines representing negative values) when $t=4021$ (d) and $t=4049$ (e). The thick coloured contours in each subfigure correspond to $\psi = 0$, $-0.25\pi U\ave(t)$, $-0.5\pi U\ave(t)$, $-0.75\pi U\ave(t)$, and $-\pi U\ave(t)$ from bottom to top.}
\label{fig:limit_cycle_Fave=0.0022}
\end{figure}

\section{Discussion}
\label{sec:discussion}
Charney and DeVore \cite{Charney1979} did not realise that the equations of the vertical vorticity (\ref{eq:psi_final})--(\ref{eq:psi_bc_final}) are ill-posed, and approximated the stream function as
\begin{equation}
\psi(x,y,t) = \sqrt{2} \psiA(t) \cos y + 2 \psiK(t) \cos 2x \sin y + 2 \psiL(t) \sin 2x \sin y,\label{eq:psi_approx}
\end{equation}
which satisfies the boundary condition (\ref{eq:psi_bc_final}) naturally. Assuming $h = h_0 \cos 2x \sin y$ and $F'=\sqrt{2}k\psiAO(\sin y - 2/\pi)$, substituting (\ref{eq:psi_approx}) into (\ref{eq:psi_final}), and using the Galerkin method, they obtained
\begin{align}
& \frac{\dif \psiA}{\dif t} = \frac{8\sqrt{2}}{3\pi} h_0 \psiL - k(\psiA - \psiAO),\label{eq:psi_A}\\
& \frac{\dif \psiK}{\dif t} = -\frac{64\sqrt{2}}{15\pi} \psiA \psiL + \frac{2}{5} \beta \psiL - k \psiK,\\
& \frac{\dif \psiL}{\dif t} = \frac{64\sqrt{2}}{15\pi} \psiA \psiK - \frac{8\sqrt{2}}{15\pi} h_0 \psiA - \frac{2}{5} \beta \psiK - k \psiL.\label{eq:psi_L}
\end{align}

Table \ref{tab:psi} lists the three equilibrium states of the equations of the vertical vorticity calculated by Charney and DeVore \cite{Charney1979} for $k = 0.01$, $\beta = 0.25$, $h_0 = 0.2$, and $\psiAO = 0.2$. The average zonal velocity is
\begin{equation}
U\ave = \frac{1}{\pi} (\psi|_{y=0} - \psi|_{y=\pi}) = \frac{2\sqrt{2}}{\pi} \psiA,\label{eq:average_velocity}
\end{equation}
and the average zonal force is
\begin{equation}
F\ave = \frac{\dif U\ave}{\dif t} + k U\ave - (hV)\ave = \frac{2\sqrt{2}}{\pi} (\frac{\dif \psiA}{\dif t} + k\psiA) - h_0 \psiL = (\frac{32}{3\pi^2}-1) h_0 \psiL + \frac{2\sqrt{2}}{\pi} k \psiAO,\label{eq:average_force}
\end{equation}
where (\ref{eq:dUdt}) and (\ref{eq:psi_A}) are used. The three equilibrium states calculated by Charney and DeVore \cite{Charney1979} have neither the same average zonal velocity nor the same average zonal force (Table \ref{tab:psi}). According to the equation (\ref{eq:psi_A}), these equilibrium states satisfy $8\sqrt{2} h_0 \psiL /(3\pi) - k(\psiA - \psiAO) = 0$, which is equivalent to
\begin{equation}
F\ave = (1-\frac{3\pi^2}{32}) k U\ave + \frac{3\pi}{8\sqrt{2}}k\psiAO,\label{eq:additional_condition}
\end{equation}
where (\ref{eq:average_velocity}) and (\ref{eq:average_force}) are used.

\begin{table}[htbp]
\ttabbox[\FBwidth]
{\caption{The equilibrium states of the equations (\ref{eq:psi_A})--(\ref{eq:psi_L}) for $k = 0.01$, $\beta = 0.25$, $h_0 = 0.2$, and $\psiAO = 0.2$ \cite{Charney1979}. $U\ave$ and $F\ave$ are calculated according to (\ref{eq:average_velocity}) and (\ref{eq:average_force}), respectively.}\label{tab:psi}}
{\begin{tabular*}{10cm}{@{\extracolsep{\fill}}@{~~}ccccc}
\hline
$\psiA$ & $\psiK$ & $\psiL$ & $U\ave$ & $F\ave$\\
\hline
 0.1535   &  0.03773    & -0.001937  & 0.1382   & 0.001769 \\
 0.1212   &  0.04358    & -0.003282  & 0.1091   & 0.001748 \\
 0.02943 & -0.03088    & -0.007104  & 0.02650 & 0.001686 \\
 \hline
\end{tabular*}}
\end{table}

Strictly speaking, the three equilibrium states obtained by Charney and DeVore \cite{Charney1979}, which are plotted in Figure \ref{fig:Fave_Uave_h0} with black squares, are the multiple solutions of the equation (\ref{eq:psi_final}) under the boundary condition (\ref{eq:psi_bc_final}) and the additional condition (\ref{eq:additional_condition}). These equilibrium states would therefore correspond to the points of intersection of the solid line that represents the branch of equilibrium states and the dashed line that represents the condition (\ref{eq:additional_condition}) in Figure \ref{fig:Fave_Uave_h0}, if there were no truncation errors. However, the additional condition (\ref{eq:additional_condition}) does not have a clear physical meaning, and is just a consequence of their severely truncated spectral expansion (\ref{eq:psi_approx}).

\section{Conclusion}
\label{sec:conclusion}
In this paper, we show that the equations of the vertical vorticity are not well-posed in the study of quasi-geostrophic barotropic flows over topography. There lacks an additional condition on the average zonal force $F\ave$, which can be given explicitly as a constant or as a function of the average zonal velocity $U\ave$, or given implicitly by prescribing $U\ave$ to be a constant.

We prove that there is an asymptotically stable equilibrium state under the condition of constant $F\ave$ or constant $U\ave$ when $F'=kC(\sin y - 2/\pi)$ and $h = C\eta \cos y$, where $C$ and $\eta>-1$ are constants. Particularly, when $h=0$, the equilibrium state $U=(F\ave + F')/k$, $V=0$ is asymptotically stable in a constant $F\ave$ problem, and the equilibrium state $U=U\ave + F'/k$, $V=0$ is asymptotically stable in a constant $U\ave$ problem, which excludes the existence of multiple equilibrium states in the case of flat topography.

When the topography is $h = h_0 \cos 2x \sin y$ and the fluctuating zonal force is
$$F'=\sqrt{2}k\psiAO(\sin y - 2/\pi),$$
 we calculate the equilibrium states with a pseudo-arclength continuation method for $k = 0.01$, $\beta = 0.25$, $\psiAO = 0.2$, and $0\leq h_0 \leq 0.2$ (Figure \ref{fig:Fave_Uave}). Multiple equilibrium states with the same $F\ave$ appear only for $h_0\geq 0.02$. Their stabilities are examined with direct numerical simulations. When $h_0=0.04$, along the $F\ave$--$U\ave$ curve, the stability of the equilibrium state changes near the extreme points of $F\ave$. However, this is not the case for $h_0=0.1$ and $h_0=0.2$, where there may be Hopf bifurcation. The emergence and evolution of limit cycles deserve further studies.

\

\noindent\textbf{Acknowledgment}.  The research was partially supported by NSFC of China (11571240, 11602148).


\bibliographystyle{elsarticle-num-names}
\biboptions{sort&compress}
\bibliography{Science_China_Math_2020_0435}





\end{document}